\documentclass{amsart}

\usepackage{mathtools}
\usepackage{microtype}
\mathtoolsset{showonlyrefs}
\usepackage{verbatim}
\newtheorem{theorem}{Theorem}
\newtheorem{proposition}[theorem]{Proposition}
\newtheorem{corollary}[theorem]{Corollary}
\newtheorem{lemma}[theorem]{Lemma}
\newtheorem{definition}[theorem]{Definition}
\newtheorem{example}[theorem]{Example}
\theoremstyle{remark}

\newcommand{\R}{\mathbb{R}}
\newcommand{\Sph}{\mathbb{S}}
\newcommand{\BP}{\mathcal{BP}}
\newcommand{\Gr}{\mathrm{Gr}}
\newcommand{\dd}{\,d}
\newcommand{\vol}[2][n]{\lvert #2\rvert_{#1}}

\title[Extensions of the Busemann--Petty problem]
{Extensions of the Busemann--Petty problem for arbitrary measures}

\author[D. Galicer]{Daniel Galicer} \address{Universidad Torcuato Di Tella. Departamento de Matem\'aticas y Estad\'istica. IMAS--CONICET. Av. Figueroa Alcorta 7350, C1428 Buenos Aires, Argentina} \email{daniel.galicer@utdt.edu}

\author[J. Haddad]{Juli\'{a}n Haddad}
\address{Departamento de An\'{a}lisis Matem\'{a}tico, Universidad de Sevilla,
Sevilla, Espa\~na}
\email{jhaddad@us.es}

\author[J. Singer]{Joaqu\'{\i}n Singer}
\address{Departamento de Matem\'{a}tica, Facultad de Ciencias Exactas y
Naturales, Universidad de Buenos Aires, (1428) Buenos Aires, Argentina}
\email{jsinger@dm.uba.ar}


\begin{document}

\begin{abstract}
The classical Busemann--Petty problem asks whether smaller central
hyperplane sections of origin-symmetric convex bodies imply smaller total
volume. Zvavitch studied the analogous question when sections and bodies are
measured by two arbitrary densities.

We refine this result in three directions: we allow central sections of
arbitrary codimension; we relax the monotonicity requirement on the radial
density ratio to a decomposition into a non-decreasing and a non-increasing
part; and we permit a distinct pair of densities for each body, one for the
sections and another for the full volume. We also obtain an isomorphic
version, in which the comparison constant is governed by the distance from an
auxiliary star body to the class of generalized $k$-intersection bodies, and
we present some examples illustrating cases not covered by previous results.
\end{abstract}

\maketitle

\section{Introduction}

Understanding how lower-dimensional information determines global geometric
quantities is a central theme in convex and geometric analysis. A prototypical
example of this phenomenon is provided by the \emph{Busemann--Petty problem}
in $\R^n$, which compares $(n-1)$-dimensional central hyperplane sections of
two origin-symmetric convex bodies and asks whether smaller sections
necessarily imply smaller $n$-dimensional volume.

We recall that a \emph{convex body} in $\R^n$ is a
compact convex set with non-empty interior. A \emph{star body} is a set
$K\subseteq\R^n$ of the form
\[
 K=\{0\}\cup\{r\theta: \theta\in\Sph^{n-1},\ 0<r\leq \rho_K(\theta)\},
\]
where the radial function $\rho_K:\Sph^{n-1}\to(0,\infty)$ is continuous.
Every convex body having the origin in its interior is a star body.
The Minkowski functional of a star body is
\[
 \lVert x\rVert_K
 =\inf\{t>0:x\in tK\}
 =\frac{|x|}{\rho_K(x/|x|)},\qquad x\neq0.
\]
It is a norm when $K$ is origin-symmetric and convex.

The classical problem posed by Busemann and Petty in \cite{BusemannPetty1956}
asks whether, for origin-symmetric convex bodies $K,L\subset\R^n$,
\[
 \vol[n-1]{K\cap\xi^\perp}
 \leq \vol[n-1]{L\cap\xi^\perp}
 \quad\text{for every }\xi\in\Sph^{n-1}
\]
implies $\vol{K}\leq\vol{L}$. The answer is affirmative for $n\leq4$ and
negative for $n\geq5$. Intersection bodies, introduced by Lutwak
\cite{Lutwak1988}, and their Fourier-analytic characterization
\cite{Koldobsky1998} are central to the solution; see \cite{Gardner1994} and
\cite{Zhang1999} for the affirmative answers in dimensions three and four,
\cite{GKS1999} for a unified analytic solution, and
\cite{GardnerBook,KoldobskyBook} for the history.

Zvavitch considered a version in which one density measures central
hyperplane sections and another measures total mass. The following form of
his theorem is convenient for comparison with our result. We recall that a
tempered distribution is \emph{positive} when it is non-negative on
non-negative Schwartz functions, and \emph{positive definite} when its Fourier
transform is positive.

\begin{theorem}[Zvavitch \cite{Zvavitch2005}]\label{thm:zvavitch}
Let $f_{n-1}$ and $f_n$ be continuous, non-negative, even functions on
$\R^n$ such that
\begin{equation}\label{eq:zvavitch-ratio}
 t\longmapsto \frac{f_{n-1}(t\theta)}{t f_n(t\theta)}
\end{equation}
is non-increasing for every $\theta\in\Sph^{n-1}$. Suppose that $K$ is an
origin-symmetric star body and that
\[
 x\longmapsto
 \lVert x\rVert_K^{-1}
 \frac{f_n(x/\lVert x\rVert_K)}
      {f_{n-1}(x/\lVert x\rVert_K)}
\]
is a positive-definite distribution. If $L$ is an origin-symmetric star body
and
\[
 \int_{K\cap\xi^\perp}f_{n-1}(x)\dd x
 \leq
 \int_{L\cap\xi^\perp}f_{n-1}(x)\dd x
 \quad\text{for every }\xi\in\Sph^{n-1},
\]
then
\[
 \int_K f_n(x)\dd x\leq\int_L f_n(x)\dd x.
\]
\end{theorem}

Theorem~\ref{thm:main} below extends Theorem~\ref{thm:zvavitch} in three ways.
First, the monotonicity hypothesis is replaced by a decomposition: the radial
quotient in \eqref{eq:zvavitch-ratio} is written as a sum of a non-decreasing
and a non-increasing term, which covers both monotone cases and, beyond them,
quotients that are monotone in neither sense. Second, the two bodies need not
share a measurement scheme --- each carries its own density for sections and
its own density for total mass, so that four measures are involved rather than
two. Third, the comparison is made for central sections of arbitrary
codimension $k$, with generalized $k$-intersection bodies in place of
intersection bodies. The conclusion is quantitative: the two total masses are
compared with a constant given by the distance from an auxiliary star body,
built from the decomposition, to the class of generalized $k$-intersection
bodies.

Both the hypothesis and the conclusion are stated on the symmetric difference
$K\mathbin\triangle L$ rather than on the bodies themselves. When a single pair
of densities is used for both bodies, the contributions of $K\cap L$ are common
to the two sides of each inequality and cancel, so that nothing is lost and the
difference form is equivalent to the usual one. When each body carries its own
pair of densities, the difference form is the natural formulation, since there
is no reason for the two schemes to agree on the overlap. It also has an
analytic benefit: $K\mathbin\triangle L$ is separated from the origin, so the
densities used for total mass are never integrated near the origin. Local
integrability there, which is required in Theorem~\ref{thm:zvavitch}, plays no
role, and densities as singular as $|x|^{-n}$ become admissible; see
Example~\ref{ex:singular}.

This general setup also provides a natural framework for isomorphic comparison
problems. In the case of Lebesgue volume, the affirmative resolution of the
slicing problem by Klartag and Lehec \cite{KlartagLehec2025} yields a
dimension-free isomorphic constant. For general measures, Koldobsky and
Zvavitch \cite{KoldobskyZvavitch2015} established a $\sqrt{n}$ estimate, which
they have recently shown to be of the sharp order
\cite{KoldobskyZvavitch2026}. Isomorphic
two-measure comparisons in terms of geometric distance to generalized
intersection bodies were subsequently studied in
\cite{KoldobskyPaourisZvavitch2022}. Our results complement these
developments: the density used on sections may differ from the density used on
total volume, and each body may be equipped with its own distinct pair of
measures.

The paper is organized as follows. Section~\ref{sec:main} contains the
background on generalized $k$-intersection bodies, the statement of the main
theorem, Theorem~\ref{thm:main}, and the isomorphic estimate that follows from
it. Section~\ref{sec:applications} collects three applications: a calibration
principle relating the density used on sections to the density used for total
mass; a family of section densities whose radial quotient attains an interior
minimum, so that neither monotonicity assumption is available; and an example
in which the density used for total mass is not locally integrable at the
origin. Section~\ref{sec:proof} contains the proof of
Theorem~\ref{thm:main}.

\section{Generalized intersection bodies and the main theorem}
\label{sec:main}

For $1\leq k<n$, let $\Gr_{n-k}$ be the Grassmannian of
$(n-k)$-dimensional subspaces of $\R^n$. The spherical Radon transform is
defined by
\[
 \mathcal R_{n-k}\varrho(H)
 =\int_{\Sph^{n-1}\cap H}\varrho(\theta)\dd\theta,
 \qquad H\in\Gr_{n-k},
\]
for $\varrho\in C(\Sph^{n-1})$.

\begin{definition}\label{def:BP}
An origin-symmetric star body $D$ belongs to $\BP_k^n$ if there is a finite
non-negative Borel measure $\eta_D$ on $\Gr_{n-k}$ such that
\begin{equation}\label{eq:BP-definition}
 \int_{\Sph^{n-1}}\rho_D(\theta)^k \varrho(\theta)\dd\theta
 =\int_{\Gr_{n-k}}\mathcal R_{n-k}\varrho(H)\dd\eta_D(H)
\end{equation}
for every even $\varrho\in C(\Sph^{n-1})$.
\end{definition}

The class $\BP_k^n$ is the class of generalized $k$-intersection bodies,
introduced by Zhang \cite{Zhang1996}; the formulation in
Definition~\ref{def:BP} is the one used by Koldobsky
\cite{Koldobsky2000}, see also \cite{KoldobskyBook}. For $k=1$ it coincides with
the usual class of intersection bodies, introduced by Lutwak
\cite{Lutwak1988} and subsequently developed by Gardner, Zhang, and others.
Koldobsky's characterization says that $D\in\BP_1^n$ if and only if
$\lVert\cdot\rVert_D^{-1}$ is a positive-definite distribution
\cite{Koldobsky1998}. Thus Definition~\ref{def:BP} recovers precisely the
analytic condition appearing in the hyperplane problem.

For an origin-symmetric star body $M$, define its distance to $\BP_k^n$ by
\begin{equation}\label{eq:distance}
 d_{\mathrm{BM}}(M,\BP_k^n)
 =\inf\bigl\{\lambda\geq1:
 \lambda^{-1}D\subseteq M\subseteq D
 \text{ for some }D\in\BP_k^n\bigr\}.
\end{equation}
Since $\BP_k^n$ is invariant under invertible linear maps, this is the usual
Banach--Mazur distance written in a normalized form.

We now state the main result. Densities indexed by $n-k$ are integrated on
$(n-k)$-dimensional subspaces, whereas densities indexed by $n$ are integrated
in $\R^n$. Throughout, if $\mu$ is a measure with density $f$, if $A\subseteq
\R^n$ is a Borel set and if $H\in\Gr_{n-k}$, we write
\[
 \mu(A\cap H)=\int_{A\cap H}f(x)\dd x,
\]
the integral being taken with respect to the $(n-k)$-dimensional Lebesgue
measure on $H$.

\begin{theorem}\label{thm:main}
Let $1\leq k<n$, let $K,L\subset\R^n$ be origin-symmetric star bodies, and
let $\mu_n,\mu_{n-k},\nu_n,\nu_{n-k}$ be four measures with non-negative even densities
$f_n,f_{n-k},g_n,g_{n-k}$, respectively. Assume that $f_{n-k}$ and
$g_{n-k}$ are continuous and that $f_n,g_n$ are integrable on
$K\mathbin\triangle L$.

On $K\mathbin\triangle L$, set
\[
 u_j(x)=
 \begin{cases}
  f_j(x),&x\in K\setminus L,\\
  g_j(x),&x\in L\setminus K,
 \end{cases}
 \qquad j\in\{n-k,n\}.
\]
Suppose that there are finite, even, measurable functions
$a,b:\R^n\setminus\{0\}\to\R$ such that $a(r\theta)$ is non-decreasing and
$b(r\theta)$ is non-increasing in $r>0$, for every
$\theta\in\Sph^{n-1}$, and
\begin{equation}\label{eq:factorization}
 u_{n-k}(x)=|x|^k u_n(x)\bigl(a(x)+b(x)\bigr)
 \quad\text{for almost every }x\in K\mathbin\triangle L.
\end{equation}
Assume further that
\begin{equation}\label{eq:M}
 c(\theta)
 =a(\rho_L(\theta)\theta)+b(\rho_K(\theta)\theta)
\end{equation}
is positive and continuous and let $M$ be the auxiliary star body with
\[
 \rho_M(\theta)=c(\theta)^{-1/k}.
\]
If
\begin{equation}\label{eq:section-hypothesis}
 \mu_{n-k}\bigl((K\setminus L)\cap H\bigr)
 \leq
 \nu_{n-k}\bigl((L\setminus K)\cap H\bigr)
 \qquad\text{for every }H\in\Gr_{n-k},
\end{equation}
then
\begin{equation}\label{eq:main-conclusion}
 \mu_n(K\setminus L)
 \leq d_{\mathrm{BM}}(M,\BP_k^n)^k\,
       \nu_n(L\setminus K).
\end{equation}
In particular, if $M\in\BP_k^n$, then
\[
 \mu_n(K\setminus L)\leq\nu_n(L\setminus K).
\]
\end{theorem}

Some important remarks are in order.

Writing the hypothesis as \eqref{eq:factorization}, rather than as a quotient,
keeps it meaningful at points where $u_n$ vanishes. As observed in the
introduction, $K$ and $L$ both contain the origin in their interiors, so
$K\mathbin\triangle L$ is separated from the origin and no local integrability
of $f_n$ or $g_n$ there is required.

The decomposition required in \eqref{eq:factorization} is available whenever
the relevant radial function has bounded variation: on the radial interval
with endpoints $\rho_K(\theta)$ and $\rho_L(\theta)$, the Jordan decomposition
writes a function of bounded variation as the sum of a non-decreasing and a
non-increasing real-valued function. The two summands need not be
non-negative. This is why Theorem~\ref{thm:main} only requires the endpoint
combination \eqref{eq:M} to be positive and continuous. Consequently, raywise
bounded variation is a convenient source of decompositions, but it must be
accompanied by this endpoint condition. In particular, global Lipschitz
regularity gives the required bounded variation on every radial interval,
while bounded variation regularity on an annulus gives it for almost every
direction.

Notice also that for $k=1$, and when the same pair of densities is assigned to
both bodies, Theorem~\ref{thm:main} directly recovers Zvavitch's result
(Theorem~\ref{thm:zvavitch}). Specifically, if the quotient in
\eqref{eq:zvavitch-ratio} is non-increasing, one takes $a=0$, so that
$c(\theta)=b(\rho_K(\theta)\theta)$ and $M$ is exactly the body appearing in
Theorem~\ref{thm:zvavitch}; if it is non-decreasing, one takes $b=0$. The
general factorization \eqref{eq:factorization} unifies these two behaviors,
allowing non-monotonic density ratios where increasing and decreasing radial
components occur simultaneously.

We point out that the quotient in \eqref{eq:zvavitch-ratio} is the reciprocal
of the one appearing in Zvavitch's original statement, where a non-decreasing
condition is required. We have chosen this normalization to make the analogy
with \eqref{eq:factorization} transparent.

Finally, the distance term gives an immediate isomorphic consequence.

\begin{corollary}\label{cor:john}
Under the hypotheses of Theorem~\ref{thm:main}, if $M$ is convex, then
\[
 \mu_n(K\setminus L)\leq n^{k/2}\nu_n(L\setminus K).
\]
In particular, for central hyperplanes,
\[
 \mu_n(K\setminus L)\leq\sqrt n\,\nu_n(L\setminus K).
\]
\end{corollary}

\begin{proof}
By John's theorem there is an origin-centered ellipsoid $E$ such that
$E\subseteq M\subseteq\sqrt n\,E$. Every ellipsoid belongs to $\BP_k^n$.
Taking $D=\sqrt n\,E$ in \eqref{eq:distance} gives
$d_{\mathrm{BM}}(M,\BP_k^n)\leq\sqrt n$.
\end{proof}

\section{Applications and examples}
\label{sec:applications}

\subsection{A calibration principle}

The first family explains a natural relation between the density used on
sections and the density used for total mass. The factor
$\lVert x\rVert_D^k$ compensates for the difference between the radial
Jacobians $r^{n-k-1}$ on an $(n-k)$-dimensional section and $r^{n-1}$ in
$\R^n$.

\begin{proposition}[Calibrated pairs of densities]\label{prop:calibration}
Let $1\leq k<n$, let $D\in\BP_k^n$, and let $\varphi,\psi$ be positive,
continuous, even functions on $\R^n$. For arbitrary origin-symmetric star
bodies $K,L$, assume that
\begin{equation}\label{eq:calibrated-sections}
 \int_{(K\setminus L)\cap H}
       \lVert x\rVert_D^k\varphi(x)\dd x
 \leq
 \int_{(L\setminus K)\cap H}
       \lVert x\rVert_D^k\psi(x)\dd x
\end{equation}
for every $H\in\Gr_{n-k}$. Then
\begin{equation}\label{eq:calibrated-total}
 \int_{K\setminus L}\varphi(x)\dd x
 \leq
 \int_{L\setminus K}\psi(x)\dd x.
\end{equation}
\end{proposition}

\begin{proof}
Use the volume densities $f_n=\varphi$ and $g_n=\psi$, and the section
densities
\[
 f_{n-k}(x)=\lVert x\rVert_D^k\varphi(x),
 \qquad
 g_{n-k}(x)=\lVert x\rVert_D^k\psi(x).
\]
Since
\[
 \frac{\lVert r\theta\rVert_D^k}{r^k}
 =\rho_D(\theta)^{-k},
\]
the factorization \eqref{eq:factorization} holds with
$a(r\theta)=\rho_D(\theta)^{-k}$, which is constant in $r$, and $b=0$. The
auxiliary body in \eqref{eq:M} is exactly $D$, so Theorem~\ref{thm:main}
applies with constant one.
\end{proof}

For example, take $k=1$, $D=B_2^n$, $\varphi(x)=e^{-|x|^2}$, and
$\psi\equiv1$. Proposition~\ref{prop:calibration} then says that
\begin{equation}\label{eq:gaussian-sections}
 \int_{(K\setminus L)\cap\xi^\perp}|x|e^{-|x|^2}\dd x
 \leq
 \int_{(L\setminus K)\cap\xi^\perp}|x|\dd x
 \quad\text{for all }\xi\in\Sph^{n-1}
\end{equation}
implies
\begin{equation}\label{eq:gaussian-example}
 \int_{K\setminus L}e^{-|x|^2}\dd x
 \leq\vol{L\setminus K}.
\end{equation}
This compares Gaussian mass on one side with Lebesgue mass on the other. It
is not an instance of Theorem~\ref{thm:zvavitch}: the two bodies are measured
by different pairs of densities. For generic non-nested bodies, patching the
two section densities across the radial interface produces a discontinuous
density, whereas Zvavitch's theorem requires one continuous section density.

\subsection{A non-monotone family}

We next describe a family of section densities that uses both extensions of
Zvavitch's theorem at once: the radial quotient is not monotone, and the two
bodies are measured by different densities.

Given $1\leq k<n$, $0<\alpha<k$ and origin-symmetric star bodies $K,L$,
define, for $x=r\theta\neq0$,
\begin{equation}\label{eq:omega}
 \omega_{K,L,\alpha}(r\theta)
 =r^k\left[
       \left(\frac{r}{\rho_L(\theta)}\right)^\alpha
       +\left(\frac{\rho_K(\theta)}{r}\right)^\alpha
     \right].
\end{equation}
Because $0<\alpha<k$, this function extends continuously to the origin by
$\omega_{K,L,\alpha}(0)=0$.

\begin{proposition}[Non-monotone section densities]
\label{prop:nonmonotone}
Let $1\leq k<n$, let $0<\alpha<k$, and let $\varphi,\psi$ be positive,
continuous, even functions on $\R^n$. For
arbitrary origin-symmetric star bodies $K,L\subset\R^n$, if
\begin{equation}\label{eq:nonmonotone-sections}
 \int_{(K\setminus L)\cap H}
     \omega_{K,L,\alpha}(x)\varphi(x)\dd x
 \leq
 \int_{(L\setminus K)\cap H}
     \omega_{K,L,\alpha}(x)\psi(x)\dd x
\end{equation}
for every $H\in\Gr_{n-k}$, then
\begin{equation}\label{eq:nonmonotone-total}
 \int_{K\setminus L}\varphi(x)\dd x
 \leq
 \int_{L\setminus K}\psi(x)\dd x.
\end{equation}
\end{proposition}

\begin{proof}
Take $f_n=\varphi$, $g_n=\psi$ and
\[
 f_{n-k}=\omega_{K,L,\alpha}\varphi,
 \qquad
 g_{n-k}=\omega_{K,L,\alpha}\psi.
\]
The factorization \eqref{eq:factorization} holds with
\[
 a(r\theta)=\left(\frac{r}{\rho_L(\theta)}\right)^\alpha,
 \qquad
 b(r\theta)=\left(\frac{\rho_K(\theta)}{r}\right)^\alpha.
\]
Here $a$ is radially non-decreasing, $b$ is radially non-increasing, and
\[
 a(\rho_L(\theta)\theta)
 +b(\rho_K(\theta)\theta)=2.
\]
Thus the auxiliary body is the Euclidean ball of radius $2^{-1/k}$, which
belongs to $\BP_k^n$. Theorem~\ref{thm:main} gives
\eqref{eq:nonmonotone-total} with constant one.
\end{proof}

This family is genuinely beyond the monotonicity assumption in
Theorem~\ref{thm:zvavitch}. Indeed, in every direction for which
$\rho_K(\theta)\neq\rho_L(\theta)$, the radial quotient
\[
 r\longmapsto
 \left(\frac{r}{\rho_L(\theta)}\right)^\alpha
 +\left(\frac{\rho_K(\theta)}{r}\right)^\alpha
\]
has its unique minimum at
$r=\sqrt{\rho_K(\theta)\rho_L(\theta)}$, strictly between the two radial
endpoints. Hence it is neither non-decreasing nor non-increasing on the
relevant interval.

In the hyperplane case, one may take $\alpha=1/2$. Then
\[
 \omega_{K,L,1/2}(r\theta)
 =\frac{r^{3/2}}{\rho_L(\theta)^{1/2}}
  +\rho_K(\theta)^{1/2}r^{1/2},
\]
which gives a particularly simple explicit example valid for all $n\geq2$.

\subsection{A singular total density}

The difference formulation also permits total densities that are not locally
integrable at the origin. The conclusion is nevertheless finite because the
symmetric difference stays away from the origin.

\begin{example}\label{ex:singular}
Let $K,L\subset\R^n$ be origin-symmetric star bodies such that the star body
$M$ with
\[
 \rho_M(\theta)=\rho_L(\theta)^{1-n}
\]
is an intersection body. If
\[
 \vol[n-1]{K\cap\xi^\perp}
 \leq\vol[n-1]{L\cap\xi^\perp}
 \qquad\text{for every }\xi\in\Sph^{n-1},
\]
then
\[
 \int_{K\setminus L}|x|^{-n}\dd x
 \leq
 \int_{L\setminus K}|x|^{-n}\dd x.
\]
\end{example}

\begin{proof}
Take $f_{n-1}=g_{n-1}=1$ and $f_n=g_n=|x|^{-n}$. Then
\[
 \frac{f_{n-1}(r\theta)}{r f_n(r\theta)}=r^{n-1}.
\]
Thus $a(x)=|x|^{n-1}$ and $b=0$, so that
$c(\theta)=a(\rho_L(\theta)\theta)=\rho_L(\theta)^{n-1}$ and the auxiliary
radial function is $\rho_L^{1-n}$. The conclusion follows from
Theorem~\ref{thm:main} with $k=1$.
\end{proof}

\section{Proof of Theorem~\ref{thm:main}}
\label{sec:proof}

Before proceeding to the proof of our main theorem, we establish a simple but
useful continuity result that will be needed in what follows.

\begin{lemma}\label{lem:continuity}
Let $f,g:\R^n\to\R$ be continuous and let $K,L$ be star bodies. Define
$Q:\Sph^{n-1}\to\R$ by
\[
 Q(\theta)=\int_{\rho_K(\theta)}^{\rho_L(\theta)}w(r\theta)\dd r,
\]
where
\[
 w(x)=
 \begin{cases}
  f(x)&\text{if }x\in K\setminus L,\\
  g(x)&\text{if }x\in L\setminus K,\\
  0&\text{otherwise.}
 \end{cases}
\]
Then $Q$ is continuous.
\end{lemma}

\begin{proof}
Let $\theta_j\in\Sph^{n-1}$ be a sequence converging to
$\theta\in\Sph^{n-1}$ as $j\to\infty$. Note that $w$ is bounded on the
compact set $K\cup L$ and vanishes outside it; let $C$ be a bound for $|w|$.
We distinguish two cases, $\rho_K(\theta)=\rho_L(\theta)$ and
$\rho_K(\theta)\neq\rho_L(\theta)$.

In the first case we have $Q(\theta)=0$, and
\[
 |Q(\theta_j)|
 \leq C\,|\rho_L(\theta_j)-\rho_K(\theta_j)|
 \longrightarrow C\,|\rho_L(\theta)-\rho_K(\theta)|=0
\]
as $j\to\infty$.

In the second case we may assume without loss of generality that
$\rho_K(\theta)<\rho_L(\theta)$; the opposite case is identical after
exchanging the roles of $f$ and $g$ and reversing the orientation of the
integral. By the continuity of $\rho_K,\rho_L$ we have
$\rho_K(\theta_j)<\rho_L(\theta_j)$ for $j$ sufficiently large. Moreover, for
any $r\in(\rho_K(\theta),\rho_L(\theta))$ we have
$\rho_K(\theta_j)<r<\rho_L(\theta_j)$ for $j$ sufficiently large. For any such
$r$ we have $w(r\theta_j)=g(r\theta_j)$, which converges to
$g(r\theta)=w(r\theta)$ as $j\to\infty$, by the continuity of $g$. Then
\[
 w(r\theta_j)\,\mathbf 1_{(\rho_K(\theta_j),\rho_L(\theta_j))}(r)
 \longrightarrow
 w(r\theta)\,\mathbf 1_{(\rho_K(\theta),\rho_L(\theta))}(r)
\]
for every $r>0$ with $r\neq\rho_K(\theta)$ and $r\neq\rho_L(\theta)$. Since
these functions are bounded by $C$ and supported in a fixed bounded interval,
dominated convergence gives
\[
 Q(\theta_j)=\int_{\rho_K(\theta_j)}^{\rho_L(\theta_j)}w(r\theta_j)\dd r
 \longrightarrow
 \int_{\rho_K(\theta)}^{\rho_L(\theta)}w(r\theta)\dd r=Q(\theta),
\]
and the lemma is proven.
\end{proof}

\begin{proof}[Proof of Theorem~\ref{thm:main}]
Define
\begin{equation}\label{eq:qtheta}
 q(\theta)=
 \int_{\rho_K(\theta)}^{\rho_L(\theta)}
 r^{n-k-1}u_{n-k}(r\theta)\dd r.
\end{equation}
By Lemma~\ref{lem:continuity}, applied to
$|x|^{n-k-1}f_{n-k}(x)$ and
$|x|^{n-k-1}g_{n-k}(x)$, the function $q$ is continuous and even. The
hypothesis
\[
 \mu_{n-k}\bigl((K\setminus L)\cap H\bigr)
 \leq
 \nu_{n-k}\bigl((L\setminus K)\cap H\bigr)
 \qquad\text{for every }H\in\Gr_{n-k}
\]
is equivalent to
\begin{equation}\label{eq:section-hypothesis2}
 0\leq
 \int_{(L\setminus K)\cap H}u_{n-k}(x)\dd x
 -\int_{(K\setminus L)\cap H}u_{n-k}(x)\dd x
 \qquad\text{for every }H\in\Gr_{n-k}.
\end{equation}
Using polar coordinates in each $H\in\Gr_{n-k}$ shows that
\eqref{eq:section-hypothesis2} is equivalent to
\[
 0\leq\mathcal R_{n-k}q(H)
 \qquad\text{for every }H\in\Gr_{n-k}.
\]

Let $D\in\BP_k^n$ and denote by $\eta_D$ its representing
measure. Integrating the last inequality and using \eqref{eq:BP-definition},
we obtain
\begin{equation}\label{eq:positive-D}
 0\leq\int_{\Gr_{n-k}}\mathcal R_{n-k}q(H)\dd\eta_D(H)
 =\int_{\Sph^{n-1}}\rho_D(\theta)^k q(\theta)\dd\theta.
\end{equation}
By \eqref{eq:factorization},
\[
 q(\theta)
 =\int_{\rho_K(\theta)}^{\rho_L(\theta)}
 r^{n-1}u_n(r\theta)\bigl(a(r\theta)+b(r\theta)\bigr)\dd r.
\]
Radial monotonicity gives, with $c$ as in \eqref{eq:M},
\begin{equation}\label{eq:pointwise-bound}
 \int_{\rho_K(\theta)}^{\rho_L(\theta)}
 r^{n-1}u_n(r\theta)\bigl(a(r\theta)+b(r\theta)\bigr)\dd r
 \leq c(\theta)
 \int_{\rho_K(\theta)}^{\rho_L(\theta)}
 r^{n-1}u_n(r\theta)\dd r.
\end{equation}
Indeed, if $\rho_K(\theta)\leq\rho_L(\theta)$, then for every
$r\in[\rho_K(\theta),\rho_L(\theta)]$ we have
$a(r\theta)\leq a(\rho_L(\theta)\theta)$ and
$b(r\theta)\leq b(\rho_K(\theta)\theta)$, so the factor $a+b$ is bounded above
by $c(\theta)$ on the interval of integration, whose orientation is positive.
If $\rho_L(\theta)\leq\rho_K(\theta)$, the same two monotonicity properties
bound $a+b$ from below by $c(\theta)$ on $[\rho_L(\theta),\rho_K(\theta)]$,
and reversing the orientation of the integral yields exactly the same
inequality.

Combining \eqref{eq:positive-D}, \eqref{eq:pointwise-bound}, and
$c=\rho_M^{-k}$ yields
\[
 0
 \leq
 \int_{\Sph^{n-1}}
 \left(\frac{\rho_D(\theta)}{\rho_M(\theta)}\right)^k
 \int_{\rho_K(\theta)}^{\rho_L(\theta)}
 r^{n-1}u_n(r\theta)\dd r\dd\theta,
\]
which, written in Cartesian coordinates, is
\begin{equation}\label{eq:before-BM}
 \int_{K\setminus L}
 \left(\frac{\rho_D(x/|x|)}{\rho_M(x/|x|)}\right)^k
 f_n(x)\dd x
 \leq
 \int_{L\setminus K}
 \left(\frac{\rho_D(x/|x|)}{\rho_M(x/|x|)}\right)^k
 g_n(x)\dd x.
\end{equation}

Now suppose that $\lambda^{-1}D\subseteq M\subseteq D$. Then
\[
 1\leq\frac{\rho_D(\theta)}{\rho_M(\theta)}\leq\lambda
 \qquad(\theta\in\Sph^{n-1}),
\]
and consequently
\[
 \mu_n(K\setminus L)
 \leq\lambda^k\nu_n(L\setminus K).
\]
Taking the infimum over all admissible $D$ and $\lambda$ proves
\eqref{eq:main-conclusion}.
\end{proof}

\bibliographystyle{abbrv}
\bibliography{biblio}

\end{document}